\newcommand{\rem}[1]{}
\documentclass{amsart}
\usepackage{amsfonts,amssymb,amsmath,amsthm,mathrsfs}
\usepackage{url}
\usepackage[dvips]{epsfig}
\newtheorem{thrm}{Theorem} 

\newtheorem{prop}[thrm]{Proposition}

\theoremstyle{definition}
\newtheorem{definition}[thrm]{Definition}
\begin{document}
\author[C.~A.~Mantica and L.~G.~Molinari]
{Carlo alberto Mantica and Luca~Guido~Molinari}
\address{L.~G.~Molinari (corresponding author): Physics Department,
Universit\`a degli Studi di Milano and I.N.F.N. sez. Milano,
Via Celoria 16, 20133 Milano, Italy -- C.~A.~Mantica: I.I.S. Lagrange, Via L. Modignani 65, 
20161, Milano, Italy.}
\email{luca.molinari@unimi.it, carloalberto.mantica@libero.it}
\subjclass[2010]{Primary 53B30, Secondary 53C25}
\keywords{Generalized Robertson-Walker space-time, Lorentzian concircular structure,
torse-forming vector, concircular vector.}
\title[]{A note on Concircular Structure space-times}
\begin{abstract}
In this note we show that Lorentzian Concircular Structure manifolds $(LCS)_n$ coincide with Generalized Robertson-Walker space-times.
\end{abstract}
\date{1 april 2018}
\maketitle
%
%
%

{\sf Generalized Robertson-Walker} (GRW) space-times were introduced in 1995 by Al{\'\i}as, Romero and S\'anchez \cite{[Ali951]} as the warped product $-1\times_{q^2} M^*$, 
where $(M^*,g^*)$ is a Riemannian submanifold. In other terms, they are Lorentzian manifolds 
characterised by a metric 
\begin{align}\label{(1.1)}
g_{ij}dx^i dx^j = -(dt)^2 + q(t)^2 g^*_{\mu\nu}(x^1,...,x^{n-1}) dx^\mu dx^\nu
\end{align}
They are interesting not only for geometry \cite{[San98],[San99],[Gut],[Rom1],[Ars14]}, but also for physics: they include relevant space-times such as Robertson-Walker, Einstein-de Sitter, static Einstein, de-Sitter, 
the Friedmann cosmological models. They are a wide generalization of space-times for cosmological models.

In 2003 A. A. Shaikh \cite{Shaikh03} introduced the notion of {\sf Lorentzian Concircular Structure} $(LCS)_n$. It is a Lorentzian manifold endowed with a unit time-like concircular vector field, i.e. $u^i u_i=-1$ and
\begin{align}\label{(1.2)}
\nabla_k u_j = \varphi (u_ku_j+g_{kj}) 
\end{align}
where $\varphi\neq 0$ is scalar function obeying
\begin{align}\label{(1.3)}
\nabla_j \varphi =\mu u_j
\end{align}
being $\mu $ a scalar function. Various authors studied the properties of $(LCS)_n$ manifolds \cite{ShBa05,ShBa06,Sh09,Hui13}.

We show that GRW and $(LCS)_n$ are the same space-times.\\

We recall few definitions that will be used in this note. The first ones are the definitions of ``torse-forming'' and
``concircular'' vector fields, by Yano:
\begin{definition}[Yano, \cite{[Yan1],[Yan2]}]
A vector field $X_j$ is named torse-forming if $\nabla_k X_j = \omega_k X_j+\varphi g_{kj}$, being $\varphi $ 
a scalar function and $\omega_k$ a non vanishing one-form.
It is named concircular if $\omega_k$ is a gradient or locally a gradient of a scalar function.
\end{definition}
\noindent
Fialkow gave a definition different from Yano's: 
\begin{definition}[Fialkow \cite{[Fia]}]
A vector field $X_j$ is named concircular if it satisfies $\nabla_k X_j =  \rho g_{kj}$,
being $\rho $ a scalar function.
\end{definition}

The following simple but deep result was recently proven:
\begin{thrm}[Bang-Yen Chen, \cite{[Ban14]}]
A $n>3$ dimensional Lorentzian manifold is a GRW space-time if and only if it admits a time-like concircular vector 
field (in the sense of Fialkow). 
\end{thrm}
It is worth noticing that for a unit time-like vector field, the torse-forming property by Yano becomes precisely
eq.\eqref{(1.2)}, with generic scalar field $\varphi $.
%
Based on Chen's theorem we proved:
\begin{prop}[Mantica and Molinari, \cite{[Sur],MM16}]
A $n>3$ dimensional Lorentzian manifold is a GRW space-time if and only if it admits a unit time-like 
torse-forming vector, \eqref{(1.2)}, that is also an eigenvector of the Ricci tensor.
\end{prop}
Now comes the equivalence: from \eqref{(1.2)} (holding either for GRW and $(LCS)_n$ space-times)
we evaluate 
$$R_{jkl}{}^m u_m = [\nabla_j,\nabla_k]u_l = (h_{kl} \nabla_j - h_{jl}\nabla_k) \varphi -\varphi^2
(u_j g_{kl}-u_k g_{jl}) $$
where $h_{kl}= u_ku_l + g_{kl}$. Contraction with $g^{jl}$ gives
\begin{align}
 R_k{}^m u_m = u_k [u^m\nabla_m \varphi +(n-1)\varphi^2 ] - (n-2) \nabla_k \varphi 
 \end{align}
If \eqref{(1.3)} holds, then $u_k$ is an eigenvector of the Ricci tensor, and we conclude that {\em a $(LCS)_n$ manifold is a GRW space-time}.\\
If $R_{km} u^m =\xi u_k$ it is $(n-2) \nabla_k \varphi = \alpha u_k$ for some scalar field $\alpha $, i.e. \eqref{(1.3)}
holds. Then we conclude that {\em a GRW space-time is a $(LCS)_n$ manifold}. 
\vfill
\end{document}